\documentclass[twoside,leqno,A4]{amsart}

\numberwithin{equation}{section}

\newcommand{\qdn}{\hspace*{-1.5mm}}
\newcommand{\qqdn}{\hspace*{-2.5mm}}

\newcommand{\+}{&\qqdn}%

\newcommand{\mb}[1]{\mathbb{#1}}

\newcommand{\ffnk}[4]{\left[\qdn\ba{#1}#3\\[2mm]#4\ea{\!\Big|\:#2}\right]}

\newcommand{\nnm}{\nonumber}
\newcommand{\be}{\begin{equation}}
\newcommand{\ee}{\end{equation}}
\newcommand{\ba}{\begin{array}}
\newcommand{\ea}{\end{array}}
\newcommand{\bmn}{\begin{eqnarray}}
\newcommand{\emn}{\end{eqnarray}}
\newcommand{\bnm}{\begin{eqnarray*}}
\newcommand{\enm}{\end{eqnarray*}}
\newcommand{\bln}{\begin{subequations}}
\newcommand{\eln}{\end{subequations}}

\newtheorem{thm}{Theorem}
\newtheorem{corl}{Corollary}
\newtheorem{exam}{Example}

\newcommand{\bbtm}[4]{\bibitem{kn:#1}{#2,}~\emph{#3,}~{#4.}}
\newcommand{\cito}[1]{\cite{kn:#1}}

\begin{document} 
\title{$\pi$ and other formulae implied by hypergeometric \\summation theorems}
\dedicatory{\textsc{\large Yong Sup Kim, Xiaoxia Wang and Arjun K. Rathie}}

\thanks{Y.S. Kim is grateful for support by Wonkwang University Research Fund (2012) and  X. Wang acknowledges support of Shanghai Leading Academic Discipline Project S30104.}
\subjclass[2000]{33C15, 33C20}
\keywords{$\pi$ formula,
          Gauss summation theorem,
          Bailey summation theorem,
          Watson summation theorem,
          extension summation theorem}


\begin{abstract}
By employing certain extended classical summation theorems, several surprising $\pi$ and other formulae are displayed.
\end{abstract}

\maketitle\thispagestyle{empty}
\markboth{Y.S. Kim, X. Wang and A.K. Rathie}
{$\pi$ and other formulae implied by hypergeometric summation theorems}

\section{\textbf{Introduction}}
In the usual notation, let $\mb{C}$ denote the set of complex numbers. For $\alpha_j\in\mb{C} (j=1, \cdots, p)$ and
$\beta_j\in\mb{C}\setminus \mb{Z}_0^-(\mb{Z}_0^-=\mb{Z}\cup\{0\}) (j=1, \cdots, q)$, the generalized hypergeometric function
$_pF_q$ with $p$ numerator parameters $\alpha_1, \cdots, \alpha_p$ and $q$ denominator parameters $\beta_1, \cdots, \beta_q$ is defined as \cito{bailey}
\bmn
{_p F_q}
\ffnk{ccc}{z}{\alpha_1,\+\cdots,\+\alpha_p}
              {\beta_1,\+\cdots,\+\beta_q}
=\sum_{k=0}^{\infty}
\frac{(\alpha_1)_k\cdots(\alpha_p)_k}
     {(\beta_1)_k\cdots(\beta_q)_k}
\cdot\frac{z^k}{k!},
\emn
($p, q \in\mb{N}_0=\mb{N}\cup{\{0\}}$; $p\leq q+1$; $p\leq q$ and $\mid z \mid<\infty$; $p=q+1$ and $\mid z \mid<1$;
$p=q+1$ and $\mathrm{Re}(\omega)>0$ ), where
\bmn
\omega=\sum_{j=1}^q\beta_j-\sum_{j=1}^p\alpha_j
\emn
and $(\lambda)_n$ is the Pochhammer symbol defined (for $\lambda\in\mb{C}$), in terms of the familiar Gamma function $\Gamma$, by
\bmn
(\lambda)_n=\frac{\Gamma(\lambda+n)}{\Gamma(\lambda)}
=\begin{cases}
1,\+n=0;\\[2mm]
\lambda(\lambda+1)\cdots(\lambda+n-1),\+n\in\mb{N}.
\end{cases}
\emn

It should be remarked here that whenever generalized hypergeometric functions reduce to quotients of gamma functions,
the results are very important from the applications point of view.

In the theory of hypergeometric and generalized hypergeometric series, classical summation theorems such as those of Gauss, Gauss second, Kummer and
Bailey for the series $_2F_1$; Watson, Dixon, Whipple and Saalsch\"{u}tz for the series $_3F_2$ and others play an important role.

Recently, good progress has been done in the direction of generalizing the above mentioned classical summation theorems.
For details, we refer the papers \cite{kn:kim-Rakha-rathie,kn:lavoie-grondin-rathie1,kn:lavoie-grondin-rathie2,kn:lavoie-grondin-rathie3,kn:Rakha-rathie,kn:vidnunas}.

On the other hand, formulas for $\pi$--series have been obtained by several mathematicians, see for examples, the papers
\cite{kn:adamchik-wagon, kn:bailey-borwein-plouffe, kn:borwein-2, kn:chan,kn:chu,kn:gourevith-guillera,kn:guillera2,kn:prudinkov dd,kn:wei,kn:weisstein, kn:zheng}.
But for the history and some introductive information on the formulae for $\pi$--series, we especially refer very interesting and
useful research papers by Bailey--Borwein \cito{bailey-borwein} and Guillera \cito{guillera}.

In our present investigation, we require the following summation theorems obtained earlier by Kim et al.  \cito{kim-Rakha-rathie}.
These summation theorems are included so that the paper may be self contained.

\textbf{Gauss summation theorem}
\bmn\label{gauss}
{_2F_1}\ffnk{cc}{1}{a, \+b}{\+c}=\frac{\Gamma(c)\Gamma(c-a-b)}{\Gamma(c-a)\Gamma(c-b)},
\emn
provided $\mathrm{Re}(c-a-b)>0$.

\textbf{Extension of Gauss summation theorem}
\bmn\label{gauss-e}
{_3F_2}\ffnk{ccc}{1}{a, \+b,\+d+1}{\+c+1,\+d}
=\frac{\Gamma(c+1)\Gamma(c-a-b)}{\Gamma(c-a+1)\Gamma(c-b+1)}
\big[(c-a-b)+\frac{a b}{d}\big],
\emn
for $\mathrm{Re}(d)>0$ and $\mathrm{Re}(c-a-b)>0$.

$\mathbf{Remark}$: In \eqref{gauss-e}, if we take $d=c$, we recover Gauss summation theorem \eqref{gauss}.

\textbf{Gauss second summation theorem}
\bmn\label{gauss2}
{_2F_1}\ffnk{c}{\frac12}{a, \:\:\:b}{\frac12(a+b+1)}
=\frac{\Gamma(\frac12)\Gamma(\frac12a+\frac12b+\frac12)}{\Gamma(\frac12a+\frac12)\Gamma(\frac12b+\frac12)}.
\emn

\textbf{Extension of Gauss second summation theorem}

\bmn\label{gauss2-e}
{_3F_2}\ffnk{cc}{\frac12}{a, \quad b, \+d+1}{\frac12(a+b+3),\+d}
\+=\+\frac{\Gamma(\frac12)\Gamma(\frac12a+\frac12b+\frac32)\Gamma(\frac12a-\frac12b-\frac12)}{\Gamma(\frac12a-\frac12b+\frac32)}\\
\+\times\+\Big[\frac{\frac12(a+b-1)-\frac{a b}{d}}{\Gamma(\frac12a+\frac12)\Gamma(\frac12b+\frac12)}+
\frac{\frac{a+b+1}{d}-2}{\Gamma(\frac12a)\Gamma(\frac12b)}\Big].\nnm
\emn

$\mathbf{Remark}$: In \eqref{gauss2-e}, if we take $d=\frac12(a+b+1)$, we recover Gauss second summation theorem \eqref{gauss2}.

\textbf{Bailey summation theorem}
\bmn\label{bailey}
{_2F_1}\ffnk{c}{\frac12}{a, \:1-a}{c}
=\frac{\Gamma(\frac12)\Gamma(\frac12c+\frac12)}{\Gamma(\frac12a+\frac12c)\Gamma(\frac12c-\frac12a+\frac12)}.
\emn

\textbf{Extension of Bailey summation theorem}

\bmn\label{bailey-e}
{_3F_2}\ffnk{ccc}{\frac12}{a, \+1-a, \+d+1}{\+c+1, \+d}
\+=\+2^{-c}\Gamma(\frac12)\Gamma(c+1)\\
\+\times\+\Big[\frac{\frac{2}{d}}{\Gamma(\frac12a+\frac12c)\Gamma(\frac12c-\frac12a+\frac12)}
+\frac{1-\frac{c}{d}}{\Gamma(\frac12a+\frac12c+\frac12)\Gamma(\frac12c-\frac12a+1)}\Big].\nnm
\emn

$\mathbf{Remark}$: In \eqref{bailey-e}, if we take $d=c$, we recover Bailey summation theorem \eqref{bailey}.

\textbf{Watson summation theorem}

\bmn\label{watson}
{_3F_2}\ffnk{c}{1}{a, \:\:b, \:\:c}{\frac12(a+b+1), 2c}
=\frac{\Gamma(\frac12)\Gamma(c+\frac12)\Gamma(\frac12a+\frac12b+\frac12)\Gamma(c-\frac12a-\frac12b+\frac12)}
{\Gamma(\frac12a+\frac12)\Gamma(\frac12b+\frac12)\Gamma(c-\frac12a+\frac12)\Gamma(c-\frac12b+\frac12)},
\emn
provided $\mathrm{Re}(2c-a-b)>-1$.

\textbf{Extension of Watson summation theorem}

\bmn\label{watson-e}
\+{_4F_3}\+\ffnk{cc}{1}{a, \quad b, \quad c \+d+1}{\frac12(a+b+1), 2c+1, \+d}\\
\+\+=\frac{2^{a+b-2}\Gamma(c+\frac12)\Gamma(\frac12a+\frac12b+\frac12)\Gamma(c-\frac12a-\frac12b+\frac12)}{\Gamma(\frac12)\Gamma(a)\Gamma(b)}\nnm\\
\+\+\times\Big[\frac{\Gamma(\frac12a)\Gamma(\frac12b)}{\Gamma(c-\frac12a+\frac12)\Gamma(c-\frac12b+\frac12)}
+(\frac{2c-d}{d})\frac{\Gamma(\frac12a+\frac12)\Gamma(\frac12b+\frac12)}{\Gamma(c-\frac12a+1)\Gamma(c-\frac12b+1)}\Big],\nnm
\emn
provided $\mathrm{Re}(d)>0$ and $\mathrm{Re}(2c-a-b)>-1$.

$\mathbf{Remark}$: In \eqref{watson-e}, if we take $d=2c$, we recover Watson summation theorem \eqref{watson}.

\section{\textbf{Main results}}
Our main results are given in the following theorems, corollaries and examples.

\subsection{Summation formulae implied by Gauss summation theorem \eqref{gauss}} \

Letting $a=\frac12+m$, $b=\frac12-m$ and $c=\frac32+m$ in \eqref{gauss}, we achieve the identity.
\begin{thm} For $m\in\mb{N}_0$, there holds the summation formulae for $\pi$
\bnm
{_2F_1}\ffnk{cc}{1}{\frac12+m,\+\frac12-m}{\+\frac32+m}=\frac{\pi}{2^{2m+1}}\frac{(\frac32)_m}{m!}.
\enm
\end{thm}

\begin{exam}[m=0]
\bnm
{_2F_1}\ffnk{cc}{1}{\frac12,\+\frac12}{\+\frac32}=\frac{\pi}{2}.
\enm
\end{exam}

\begin{exam}[m=1]\label{wei2}
\bnm
{_2F_1}\ffnk{cc}{1}{\frac32,\+-\frac12}{\+\frac52}=\frac{3\pi}{16}.
\enm
\end{exam}

\begin{exam}[m=2]\label{wei3}
\bnm
{_2F_1}\ffnk{cc}{1}{\frac52,\+-\frac32}{\+\frac72}=\frac{15\pi}{256}.
\enm
\end{exam}

\begin{exam}[m=3]\label{wei4}
\bnm
{_2F_1}\ffnk{cc}{1}{\frac72,\+-\frac52}{\+\frac92}=\frac{35\pi}{2048}.
\enm
\end{exam}

\begin{exam}[m=4]
\bnm
{_2F_1}\ffnk{cc}{1}{\frac92,\+-\frac72}{\+\frac{11}{2}}=\frac{315\pi}{65536}.
\enm
\end{exam}
$\mathbf{Remark}$: Results given in Examples \ref{wei2}, \ref{wei3} and \ref{wei4} are obtained recently by Wei et al.  \cito{wei} by
employing Whipple summation theorem.

\subsection{\textbf{Summation formulae implied by extension of Gauss summation theorem \eqref{gauss-e}}} \

Letting $a=\frac12+m$, $b=\frac12-m$ and $c=\frac32+m$ in \eqref{gauss-e}, we achieve the following identity.
\begin{thm} For $m\in\mb{N}_0$ and $\mathrm{Re}(d)>0$, there holds the summation formulae for $\pi$
\bnm
{_3F_2}\ffnk{ccc}{1}{\frac12+m,\+\frac12-m,\+d+1}{\+\frac52+m,\+d}=(1+\frac{1-2m}{2d})\frac{3\pi}{2^{2m+3}m!}(\frac52)_m.
\enm
\end{thm}

\begin{corl} [m=0]
\bnm
{_3F_2}\ffnk{ccc}{1}{\frac12,\+\frac12, \+d+1}{\+\frac52,\+d}=\frac{3\pi}{8}(1+\frac{1}{2d}).
\enm
\end{corl}

\begin{exam}[d=1]
\bnm
{_3F_2}\ffnk{ccc}{1}{\frac12,\+\frac12, \+2}{\+\frac52,\+1}=\frac{9\pi}{16}.
\enm
\end{exam}

\begin{exam}[d=2]
\bnm
{_3F_2}\ffnk{ccc}{1}{\frac12,\+\frac12, \+3}{\+\frac52,\+2}=\frac{15\pi}{32}.
\enm
\end{exam}

\begin{exam}[d=3]
\bnm
{_3F_2}\ffnk{ccc}{1}{\frac12,\+\frac12, \+4}{\+\frac52,\+3}=\frac{7\pi}{16}.
\enm
\end{exam}

\begin{corl} [m=1]
\bnm
{_3F_2}\ffnk{ccc}{1}{\frac32,\+-\frac12, \+d+1}{\+\frac72,\+d}=\frac{15\pi}{64}(1-\frac{1}{2d}).
\enm
\end{corl}

\begin{exam}[d=1]
\bnm
{_3F_2}\ffnk{ccc}{1}{\frac32,\+-\frac12, \+2}{\+\frac72,\+1}=\frac{15\pi}{128}.
\enm
\end{exam}

\begin{exam}[d=2]
\bnm
{_3F_2}\ffnk{ccc}{1}{\frac32,\+-\frac12, \+3}{\+\frac72,\+2}=\frac{45\pi}{256}.
\enm
\end{exam}

\begin{exam}[d=3]
\bnm
{_3F_2}\ffnk{ccc}{1}{\frac32,\+-\frac12, \+4}{\+\frac92,\+3}=\frac{75\pi}{384}.
\enm
\end{exam}

\begin{corl} [m=2]
\bnm
{_3F_2}\ffnk{ccc}{1}{\frac52,\+-\frac32, \+d+1}{\+\frac92,\+d}=\frac{105\pi}{1024}(1-\frac{3}{2d}).
\enm
\end{corl}

\begin{exam}[d=2]
\bnm
{_3F_2}\ffnk{ccc}{1}{\frac52,\+-\frac32, \+3}{\+\frac92,\+2}=\frac{105\pi}{4096}.
\enm
\end{exam}

\begin{exam}[d=3]
\bnm
{_3F_2}\ffnk{ccc}{1}{\frac52,\+-\frac32, \+4}{\+\frac92,\+3}=\frac{105\pi}{2048}.
\enm
\end{exam}

\begin{exam}[d=4]
\bnm
{_3F_2}\ffnk{ccc}{1}{\frac52,\+-\frac32, \+5}{\+\frac92,\+4}=\frac{525\pi}{8192}.
\enm
\end{exam}

\subsection{\textbf{Summation formulae implied by extension of Gauss second summation theorem \eqref{gauss2-e}}} \

Letting $a=1+2m$ and $b=1+2n$ in \eqref{gauss2-e}, we get the following identity.

\begin{thm}
For $m, n \in\mb{N}_0$ and $\mathrm{Re}(d)>0$, there holds the summation formula
\bnm
{_3F_2}\ffnk{ccc}{\frac12}{1+2m,\+1+2n,\+d+1}{\+m+n+\frac52,\+d}
=\frac{2\pi(\frac32)_{m+n+1}}{(2m-2n+1)(2m-2n-1)}\\
\times \Big[\frac{(m+n+\frac12)-\frac{1}{d}(2m+1)(2n+1)}{m!\:n!}+\frac{\frac{1}{d}(2m+2n+3)-2}{\pi(\frac12)_m(\frac12)_n}\Big].
\enm
\end{thm}

\begin{corl}[m=n=0]
\bnm
{_3F_2}\ffnk{ccc}{\frac12}{1,\+1,\+d+1}{\+\frac52,\+d}=3\pi(\frac{1}{d}-\frac12)-3(\frac{3}{d}-2).
\enm
\end{corl}

\begin{exam}[d=2]
\bnm
{_3F_2}\ffnk{ccc}{\frac12}{1,\+1,\+3}{\+\frac52,\+2}=\frac32.
\enm
\end{exam}

\begin{exam}[d=3/2]\label{wei16}
\bnm
{_2F_1}\ffnk{cc}{\frac12}{1,\+1}{\+\frac32}=\frac{\pi}{2}.
\enm
\end{exam}

\begin{corl}[m=1, n=0]
\bnm
{_3F_2}\ffnk{ccc}{\frac12}{3,\+1,\+d+1}{\+\frac72,\+d}=\frac{5\pi}{2}\big[3(\frac12-\frac{1}{d})+\frac{2}{\pi}(\frac{5}{d}-2)\big].
\enm
\end{corl}

\begin{exam}[d=2]
\bnm
{_3F_2}\ffnk{ccc}{\frac12}{3,\+1,\+3}{\+\frac72,\+2}=\frac52.
\enm
\end{exam}

\begin{exam}[d=5/2]\label{wei18}
\bnm
{_2F_1}\ffnk{cc}{\frac12}{3,\+1}{\+\frac52}=\frac{3\pi}{4}.
\enm
\end{exam}

\begin{corl}[m=1, n=1]
\bnm
{_3F_2}\ffnk{ccc}{\frac12}{3,\+3,\+d+1}{\+\frac92,\+d}=\frac{105\pi}{4}(\frac{9}{d}-\frac{5}{2})-105(\frac{7}{d}-2).
\enm
\end{corl}

\begin{exam}[d=7/2]
\bnm
{_2F_1}\ffnk{cc}{\frac12}{3,\+3}{\+\frac72}=\frac{15\pi}{8}.
\enm
\end{exam}

\begin{corl}[m=2, n=0]
\bnm
{_3F_2}\ffnk{ccc}{\frac12}{5,\+1,\+d+1}{\+\frac92,\+d}=\frac{7\pi}{4}\Big[\frac52(\frac{1}{2}-\frac{1}{d})+\frac{4}{3\pi}(\frac{7}{d}-2) \Big].
\enm
\end{corl}

\begin{exam}[d=2]
\bnm
{_3F_2}\ffnk{ccc}{\frac12}{5,\+1,\+3}{\+\frac92,\+2}=\frac72.
\enm
\end{exam}

\begin{exam}[d=7/2]\label{wei21}
\bnm
{_2F_1}\ffnk{cc}{\frac12}{5,\+1}{\+\frac72}=\frac{15\pi}{16}.
\enm
\end{exam}

$\mathbf{Remark}$: The results given in Examples \ref{wei16}, \ref{wei18} and \ref{wei21} are also obtained recently by Wei et al. \cito{wei}.

\subsection{\textbf{Summation formulae implied by extension of Bailey summation theorem \eqref{bailey-e}}} \

Letting $a=\frac12+m$ and $c=\frac32+m+2n$ in \eqref{bailey-e}, we achieve the following identity.
\begin{thm}
For $m, n \in\mb{N}_0$ and $\mathrm{Re}(d)>0$, there holds the summation formula
\bnm
{_3F_2}\ffnk{ccc}{\frac12}{\frac12+m,\+\frac12-m,\+d+1}{\+m+2n+\frac52,\+d}
=\frac{\pi}{2^{m+2n+\frac52}}(\frac32)_{m+2n+1}\\
\times \Big[\frac{\frac{2}{d}}{(m+n)!\:n!}+\frac{1-\frac{3/2+m+2n}{d}}{\pi(\frac12)_{m+n+1}(\frac12)_{n+1}}\Big].
\enm
\end{thm}

\begin{corl}[m=n=0]
\bnm
{_3F_2}\ffnk{ccc}{\frac12}{\frac12,\+\frac12,\+d+1}{\+\frac52,\+d}=\frac{3\pi\sqrt{2}}{8d}+\frac{3\sqrt{2}}{4}(1-\frac{3}{2d}).
\enm
\end{corl}

\begin{exam}[d=3/2]
\bnm
{_2F_1}\ffnk{ccc}{\frac12}{\frac12,\+\frac12}{\+\frac32}=\frac{\pi\sqrt{2}}{4}.
\enm
\end{exam}

\begin{corl}[m=1, n=0]
\bnm
{_3F_2}\ffnk{ccc}{\frac12}{\frac32,\+-\frac12,\+d+1}{\+\frac72,\+d}=\frac{15\pi\sqrt{2}}{32d}+\frac{5\sqrt{2}}{8}(1-\frac{5}{2d}).
\enm
\end{corl}

\begin{exam}[d=5/2]
\bnm
{_2F_1}\ffnk{ccc}{\frac12}{\frac32,\+-\frac12}{\+\frac52}=\frac{3\pi\sqrt{2}}{16}.
\enm
\end{exam}

\begin{corl}[m=2, n=0]
\bnm
{_3F_2}\ffnk{ccc}{\frac12}{\frac52,\+-\frac32,\+d+1}{\+\frac92,\+d}=\frac{105\pi\sqrt{2}}{256d}+\frac{7\sqrt{2}}{16}(1-\frac{7}{2d}).
\enm
\end{corl}

\begin{exam}[d=7/2]
\bnm
{_2F_1}\ffnk{ccc}{\frac12}{\frac52,\+-\frac32}{\+\frac72}=\frac{15\pi\sqrt{2}}{128}.
\enm
\end{exam}

\begin{corl}[m=3, n=0]
\bnm
{_3F_2}\ffnk{ccc}{\frac12}{\frac72,\+-\frac52,\+d+1}{\+\frac{11}{2},\+d}=\frac{315\pi\sqrt{2}}{1024d}+\frac{9\sqrt{2}}{32}(1-\frac{9}{2d}).
\enm
\end{corl}

\begin{exam}[d=9/2]
\bnm
{_2F_1}\ffnk{ccc}{\frac12}{\frac72,\+-\frac52}{\+\frac92}=\frac{35\pi\sqrt{2}}{512}.
\enm
\end{exam}

\subsection{\textbf{Summation formulae implied by extension of Watson summation theorem \eqref{watson-e}}} \

Letting $a=1+2m$, $b=1+2n$ and $c=1+m+n+s$ in \eqref{watson-e}, we obtain the following identity.
\begin{thm}
For $m, n, s \in\mb{N}_0$ and $\mathrm{Re}(d)>0$, there holds the summation formula
\bnm
\+\+{_4F_3}\ffnk{cccc}{1}{1+2m,\+1+2n,\+1+m+n+s,\+d+1}{\+m+n+\frac32,\+3+2m+2n+2s,\+d}\\
\+\+=\frac{\pi}{4}\frac{(\frac32)_{m+n+s}(\frac32)_{m+n}(\frac12)_{s}}{(\frac12)_{m}(\frac12)_{n}m!\:\:n!}
 \Big[\pi \frac{(\frac12)_{m}(\frac12)_{n}}{(m+s)!(n+s)!}+\frac{1}{\pi}\frac{(\frac{2+2m+2n+2s}{d}-1)m!\:n!}{(\frac12)_{m+s+1}(\frac12)_{n+s+1}}\Big].
\enm
\end{thm}

\begin{corl}[m=n=s=0]
\bnm
{_4F_3}\ffnk{cccc}{1}{1,\+1,\+1,\+d+1}{\+\frac32,\+\:3,\+d}=\frac{{\pi}^2}{4}+(\frac{2}{d}-1).
\enm
\end{corl}

\begin{exam}[d=1]
\bnm
{_3F_2}\ffnk{ccc}{1}{1,\+1,\+2}{\+\frac32,\+3}=\frac{{\pi}^2}{4}+1.
\enm
\end{exam}

\begin{exam}[d=2]
\bnm
{_3F_2}\ffnk{ccc}{1}{1,\+1,\+1}{\+\frac32,\+2}=\frac{{\pi}^2}{4}.
\enm
\end{exam}

\begin{exam}[d=3]
\bnm
{_4F_3}\ffnk{cccc}{1}{1,\+1,\+1,\+4}{\+\frac32,\+3,\+3}=\frac{{\pi}^2}{4}-\frac13.
\enm
\end{exam}

\begin{corl}[m=1, n=s=0]
\bnm
{_4F_3}\ffnk{cccc}{1}{3,\+1,\+2,\+d+1}{\+\frac52,\+\: 5,\+d}=\frac{9{\pi}^2}{16}+3(\frac{4}{d}-1).
\enm
\end{corl}

\begin{exam}[d=1]
\bnm
{_3F_2}\ffnk{ccc}{1}{3,\+2,\+2}{\+\frac52,\+\: 5}=\frac{9{\pi}^2}{16}+9.
\enm
\end{exam}

\begin{exam}[d=2]
\bnm
{_3F_2}\ffnk{ccc}{1}{3,\+1,\+3}{\+\frac52,\+5}=\frac{9{\pi}^2}{16}+3.
\enm
\end{exam}

\begin{exam}[d=3]
\bnm
{_3F_2}\ffnk{ccc}{1}{1,\+2,\+4}{\+\frac52,\+5}=\frac{9{\pi}^2}{16}+1.
\enm
\end{exam}

\begin{exam}[d=4]
\bnm
{_3F_2}\ffnk{ccc}{1}{3,\+1,\+2}{\+\frac52,\+4}=\frac{9{\pi}^2}{16}.
\enm
\end{exam}

\begin{exam}[d=5]
\bnm
{_4F_3}\ffnk{cccc}{1}{3,\+1,\+2,\+6}{\+\frac52,\+5,\+5}=\frac{9{\pi}^2}{16}-\frac35.
\enm
\end{exam}

\begin{corl}[m=2, n=s=0]
\bnm
{_4F_3}\ffnk{cccc}{1}{5,\+1,\+3,\+d+1}{\+\frac72,\+\:7,\+d}=\frac{225{\pi}^2}{256}+5(\frac{6}{d}-1).
\enm
\end{corl}

\begin{exam}[d=1]
\bnm
{_3F_2}\ffnk{ccc}{1}{5,\+3,\+2}{\+\frac72,\+7}=\frac{225{\pi}^2}{256}+25.
\enm
\end{exam}

\begin{exam}[d=2]
\bnm
{_4F_3}\ffnk{cccc}{1}{5,\+1,\+3,\+3}{\+\frac72,\+7,\+2}=\frac{225{\pi}^2}{256}+10.
\enm
\end{exam}

\begin{exam}[d=3]
\bnm
{_3F_2}\ffnk{ccc}{1}{5,\+1,\+4}{\+\frac72,\+7}=\frac{225{\pi}^2}{256}+5.
\enm
\end{exam}

\begin{exam}[d=4]
\bnm
{_4F_3}\ffnk{cccc}{1}{5,\+1,\+3,\+5}{\+\frac72,\+7,\+4}=\frac{225{\pi}^2}{256}+\frac52.
\enm
\end{exam}

\begin{exam}[d=5]
\bnm
{_3F_2}\ffnk{ccc}{1}{1,\+3,\+6}{\+\frac72,\+7}=\frac{225{\pi}^2}{256}+1.
\enm
\end{exam}

\begin{exam}[d=6]
\bnm
{_3F_2}\ffnk{ccc}{1}{5,\+1,\+3}{\+\frac72,\+6}=\frac{225{\pi}^2}{256}.
\enm
\end{exam}

\begin{exam}[d=7]
\bnm
{_4F_3}\ffnk{cccc}{1}{5,\+1,\+3,\+8}{\+\frac72,\+7,\+7}=\frac{225{\pi}^2}{256}-\frac57.
\enm
\end{exam}



Yong Sup Kim\\
Department of Mathematics Education\\
Wonkwang University\\
Iksan 570-749, Korea\\
\emph{E-mail address:} \textbf{yspkim@wonkwang.ac.kr}\\[0.5mm]

Xiaoxia Wang\\
Department of Mathematics\\
Shanghai University\\
Shanghai, 200444, P. R. China \\
\emph{E-mail address}: \textbf{xiaoxiawang@shu.edu.cn}\\
Corresponding author\\[0.5mm]

Arjun K. Rathie\\
Department of Mathematics\\
School of Mathematical and Physical Sciences, Central University of Kerala\\
Riverside Transit Campus, Padennakkad P..O. Nileshwar Kasaragod-671 328, Kerala State, India\\
\emph{E-mail address}: \textbf{akrathie@rediffmail.com}
\end{document}